\documentclass{amsart}

\title{Shelling totally nonnegative flag varieties}

\author[Lauren K. Williams]{Lauren K. Williams}

\address{Department of Mathematics, University of California, Berkeley, CA 94720}

\usepackage{pstricks, pst-node, amssymb}

\psset{unit=1pt, arrowsize=4pt, linewidth=.7pt}
\psset{linecolor=blue}
\newgray{grayish}{.90}
\newrgbcolor{embgreen}{0 .5 0}
\def\Le{\hbox{\rotatedown{$\Gamma$}}}
\def\vblack(#1, #2)#3{\cnode*[linecolor=black](#1, #2){3}{#3}}
\def\vwhite(#1,#2)#3{\cnode[linecolor=black,fillcolor=white,fillstyle=solid](#1,#2){3}{#3}}
\countdef\x=23
\countdef\y=24
\countdef\z=25
\countdef\t=26

\def\tbox(#1,#2)#3{
\x=#1 \y=#2 
\multiply\x by 12 
\multiply\y by 12 
\z=\x \t=\y
\advance\z by 12 
\advance\t by 12 
\psline(\x,\y)(\x,\t)(\z,\t)(\z,\y)(\x,\y)
\advance\x by 6
\advance\y by 6 
\rput(\x,\y){{\bf #3}}}

\font\co=lcircle10

\def\jr{\rotatedown{\smash{\raise2pt\hbox{\co \rlap{\rlap{\char'005} \char'007}}
               \raise6pt\hbox{\rlap{\vrule height6.5pt}}
                \raise2pt\hbox{\rlap{\hskip4pt \vrule
          height0.4pt depth0pt
                width7.7pt}}}}}
\def\textcross{\ \smash{\lower4pt\hbox{\rlap{\hskip4.15pt\vrule height14pt}}
                \raise2.8pt\hbox{\rlap{\hskip-3pt \vrule height.4pt depth0pt
		                width14.7pt}}}\hskip12.7pt}
				
\def\textelbow{\ \hskip.1pt\smash{\raise2.75pt%
                \hbox{\co \hskip 4.15pt\rlap{\rlap{\char'004} \char'006}
                \lower6.8pt\rlap{\vrule height3.5pt}
                \raise3.6pt\rlap{\vrule height3.5pt}}
                \raise2.8pt\hbox{%
                  \rlap{\hskip-7.15pt \vrule height.4pt depth0pt 
width3.5pt}%
                  \rlap{\hskip4.05pt \vrule height.4pt depth0pt 
width3.5pt}}}
                \hskip8.7pt}

\usepackage{amsmath, amsthm, amssymb, amsbsy}
\usepackage{amsfonts, latexsym, stmaryrd, amscd, xy}
\usepackage[mathscr]{eucal}
\usepackage{epsfig}
\newtheorem{theoremB}{Theorem}
\newtheorem{theorem}{Theorem}[section]

\newtheorem{proposition}[theorem]{Proposition}
\newtheorem{lemma}[theorem]{Lemma}

\newtheorem{example}[theorem]{Example}
\newtheorem{corollary}[theorem]{Corollary}
\newtheorem{conjecture}[theorem]{Conjecture}
\newtheorem{remark}[theorem]{Remark}
\newtheorem{deftheorem}[theorem]{Definition/Theorem}
\newtheorem{definition}[theorem]{Definition}

\usepackage{amsfonts}
\usepackage{xy}
\usepackage{amssymb}
\usepackage{amsmath}

\newcommand{\R}{\mathbb R}

\newcommand{\N}{\mathbb N}
\newcommand{\C}{\mathbb C}
\newcommand{\PP}{\mathcal{P}}

\newcommand{\K}{\mathcal{K}}
\newcommand{\F}{\mathcal{F}}

\newcommand{\LL}{\blacktriangleleft}
\newcommand{\RR}{\blacktriangleright}

\newcommand{\B}{\mathcal{B}}
\newcommand{\I}{\mathcal{I}}
\newcommand{\E}{\mathcal{E}}
\newcommand{\Q}{\mathcal{Q}}

\DeclareMathOperator{\rank}{rank}

\DeclareMathOperator{\Hom}{Hom}

\DeclareMathOperator{\Sl}{SL}

\newcommand{\thmrefer}[1]{\renewcommand\thetheorem
  {\protect\ref{#1}}\addtocounter{theorem}{-1}}

\xyoption{all}
\CompileMatrices

\begin{document}

\begin{abstract}
In this paper we study the partially ordered set $\Q^J$ of 
cells in Rietsch's \cite{Rietsch2} cell decomposition of the 
totally nonnegative part of an arbitrary flag variety $\PP^J_{\geq 0}$.
Our goal is to understand the geometry of $\PP^J_{\geq 0}$: 
Lusztig \cite{Lusztig2}
has proved that this space is contractible, but it is unknown 
whether the closure of each cell is contractible, and 
whether $\PP^J_{\geq 0}$ is homeomorphic to a ball.  
The order complex $\Vert \Q^J \Vert$ is a simplicial complex which
can be thought of as a combinatorial approximation of $\PP^J_{\geq 0}$.
Using combinatorial tools such as Bjorner's EL-labellings \cite{B1}
and Dyer's reflection orders \cite{Dyer},
we prove that $\Q^J$ is graded, thin and EL-shellable.  
As a corollary, we deduce that $\Q^J$ is Eulerian and that the Euler 
characteristic of the closure of each cell is $1$.  Additionally, 
our results imply that
$\Vert \Q^J \Vert$ is homeomorphic to a ball, and moreover, that 
$\Q^J$ is the face poset of some regular CW complex homeomorphic to a ball.
\end{abstract}

\maketitle

\section{Introduction}
The classical theory of total positivity concerns matrices 
in which all minors are nonnegative.  In the past decade 
Lusztig has extended 
this subject \cite{Lusztig1,
Lusztig2, Lusztig3}
by introducing the totally nonnegative 
variety $G_{\geq 0}$ in an arbitrary reductive group $G$ as well as 
the totally nonnegative part $B_{\geq 0}$ of a real flag variety $B$, 
which he called a ``remarkable polyhedral subspace" \cite{Lusztig3}.
Rietsch has constructed a cell decomposition \cite{Rietsch1} for 
the totally nonnegative part of an arbitrary flag variety $\PP^J_{\geq 0}$,
and has described the order relation
for closures of cells.  The partially ordered set (poset)
of cells, which we denote by 
$\Q^J$, is intimately connected to the Bruhat order of the corresponding
Weyl group $W$: for example, when $\PP^J$ is the complete flag variety,
$\Q^J$ is the interval poset of the Bruhat order, and when $\PP^J$
is the type A Grassmannian, $\Q^J$ is Postnikov's {\it cyclic Bruhat order}.

Lusztig \cite{Lusztig1} has proved that the totally nonnegative part of the
(full) flag variety is contractible, which implies the same result 
for any partial flag variety \cite{Rietsch3}.  However, it is 
unknown whether these spaces are homeomorphic to balls.  In addition, 
the topology of the individual cells 
is not well understood: for example, it is unknown
whether the closure of a cell is contractible.

The goal of this paper is to apply combinatorial methods to $\Q^J$
in order to better understand
the geometry of $\PP^J_{\geq 0}$.  Indeed, the past thirty years
have seen a wealth of literature 
designed to 
facilitate the interplay between combinatorics and geometry 
(for references related to 
the very important property of {\it shellability}, see
\cite{DK}, \cite{B1}, \cite{BW1}, \cite{BW2}, \cite{B2}).
In particular,
in a 1984 paper \cite{B2}, Bjorner recognized that regular CW complexes
are combinatorial objects in the following sense:
if $Q$ is the poset of closed cells in a regular CW decomposition of 
a space $P$, then the order complex (or nerve) $\Vert Q \Vert$
is homeomorphic to $P$.  Furthermore, he gave criteria \cite{B2} 
for recognizing
when a poset is the face poset of a regular CW complex: for example,
if a poset is {\it thin} and {\it shellable} then it is the face
poset of some regular CW complex homeomorphic to a sphere.

This theory has been applied to intervals in the Bruhat order of a Coxeter
group.  Bjorner and Wachs \cite{BW1} proved that such intervals are 
thin and {\it CL-shellable}; Dyer later used {\it reflection orders} to 
prove Bjorner's conjecture
that such intervals are actually {\it EL-shellable}.  Bjorner \cite{B2} 
pointed out that such intervals are therefore face posets of regular
CW complexes homeomorphic to spheres and asked for a natural geometric
construction of such a CW complex.  Sixteen years later, Fomin and Shapiro
\cite{Fomin}
had the insight that such a construction might come from total positivity,
and conjectured that links in the Bruhat decomposition of the totally
nonnegative part of the unipotent radical of a semisimple group $G$
have the desired properties.  Although they could not prove regularity
of these stratified spaces, they proved contractibility of both the
closed and open strata in type A.

In this paper we begin a line of research that in many ways parallels the 
story of the Bruhat order.
Our main results are the following.
\begin{theorem}\label{Theorem1}
$\Q^J$ is thin.
\end{theorem}
\begin{theorem}\label{Theorem2}
$\Q^J$ is EL-shellable. 
\end{theorem}
As a corollary, we deduce that $\Q^J$ is {\it Eulerian} and that the
Euler characteristic of the closure of each cell is $1$.
Additionally, it 
follows from our main
results that $\Vert \Q^J \Vert$ is homeomorphic to a ball, and that
moreover,
$\Q^J$ 
is the face poset of {\it some} regular CW complex homeomorphic to a ball.

The structure of this paper is as follows.
We begin by reviewing tools 
from the topology of partially ordered sets, 
especially Bjorner's notion 
of EL-shellability \cite{B1}.  In 
Section \ref{Bruhat}  we review 
properties of the 
Bruhat order and define the reflection orders of Dyer \cite{Dyer}.
In Section \ref{Poset} we introduce the relevant notions from 
total positivity, and Rietsch's \cite{Rietsch1, Rietsch2} 
description 
of the poset $\Q^J$  of cells in 
her cell decomposition of the totally nonnegative part of a flag variety
$\PP^J_{\geq 0}$.
In Sections \ref{thin} and \ref{EL}
we prove our main results: that the poset $\Q^J$ is thin and 
EL-shellable.  We then use Bjorner's theorem
\cite{B2} to conclude 
that $\Q^J$ is the face poset
of a 
regular CW complex homeomorphic to a ball.
We conjecture that in fact the totally nonnegative part of an arbitrary 
flag variety is a regular CW complex homeomorphic to a ball.

Additionally, we include 
an appendix which
gives a detailed explanation of some of the 
combinatorics of $\Q^J$ in the type A Grassmannian case, relating our
present work to Postnikov's \cite{Postnikov}.

\textsc{Acknowledgements:}
I am very grateful to the following people: 
to Michelle Wachs, whose lecture course at the Park City Math Institute
inspired me to start working on this subject in earnest;
to Konni Rietsch, for replying to my emails with Latexed proofs, and 
whose work on $\PP^J_{\geq 0}$ made my work possible;
to Anders Bjorner, whose work and encouragement have been very inspiring
to me; to Sergey Fomin, for editorial and philosophical advice; to 
Mark Haiman, who told me about reflection orders; and to   
Richard Stanley and 
Alex Postnikov, for their encouragement and advice.
Finally, I would like to acknowledge the hospitality of the
Mittag-Leffler Institute
in Stockholm and of Toscanini's in Cambridge, 
where some of this work was carried out.

\section{Background on poset topology}\label{PosetTopology}

Poset topology is the study of combinatorial properties of a partially
ordered set, or {\it poset}, which reflect the topology of an 
associated simplicial complex.  In this section we will
review some of the basic definitions and results of poset topology.

Let $P$ be a poset
with order relation $<$.  We will use the symbol 
$\lessdot$ to denote the covering relation in the poset: 
$x \lessdot y$ means that $x < y$ and there is no $z$ such 
that $x < z < y$.  Additionally, if $x < y$ then 
$[x,y]$ denotes the {\it interval} from $x$ to $y$; that is, the 
set $\{z \in P \ \vert \ x \leq z \leq y \}$.

The natural geometric object that one associates to a poset $P$ is
the realization of its {\it order complex} (or {\it nerve}).
The order complex $\Delta(P)$ is defined to be the simplicial
complex whose vertices are the elements of $P$ and whose simplices are
the chains $x_0 < x_1 < \dots < x_k$ in $P$.  We will denote
the realization of $\Delta(P)$ by $\Vert \Delta(P) \Vert$, or just
$\Vert P \Vert$, and if we say that $P$ possesses some topological 
property, we will mean that $\Vert P \Vert$ possesses that property.
As we will see later, it is particularly useful
to study $\Vert P \Vert$ when $P$ is the face poset of 
a CW complex.

Recall that a poset is said to be 
{\it bounded} if it has a least element and a greatest element.  These
elements are denoted $\hat{0}$ and $\hat{1}$, respectively.  A finite 
poset is said to be {\it pure} if all maximal
chains have the same length, and {\it graded}, if in addition, it is finite
and bounded.  Any element $x$ of a graded poset $P$ has a well-defined
{\it rank} $\rho(x)$ equal to the common length of all
unrefinable chains from $\hat{0}$ to $x$ in $P$.   

We now introduce the property of being {\it thin}.

\begin{definition}
A poset $P$ is called {\it thin} if every interval of length $2$ 
is a {\it diamond},
i.e. if for any $p < q$ such that $\rank(q)-\rank(p) = 2$,
there are exactly two elements in the open interval $(p,q)$.
\end{definition}

Another useful property is that of {\it shellability}.
This property is desirable to have because it implies, for 
example, that $\Delta$ is Cohen-Macaulay.

\begin{definition}
A pure finite simplicial complex $\Delta$ is said to be {\it shellable}
if its maximal faces can be ordered $F_1, F_2, \dots , F_n$ in such 
a way that 
$F_k \cap (\cup_{i=1}^{k-1} F_i)$ is a nonempty union of maximal proper
faces of $F_k$ for $k=2, 3, \dots , n$.
\end{definition}

One technique that can be used to prove that an order complex 
$\Delta(P)$ is shellable is the notion of 
{\it lexicographic shellability}, or 
{EL-shellability}, which was first introduced by Bjorner \cite{B1}.
Let $P$ be a graded poset, and let 
$\E(P)$ be the set of edges of the Hasse diagram of $P$, i.e. 
$\E(P) = \{(x,y) \in P \times P\ \vline \ x \gtrdot y\}$.  
An {\it edge labeling} of $P$ is a map 
$\lambda : \E(P) \to \Lambda$ where $\Lambda$ is some poset 
(usually the integers).  Given an edge labeling $\lambda$, 
each maximal chain $c = (x_0 \gtrdot x_1 \gtrdot \dots \gtrdot x_k)$ of 
length $k$ can be associated with a $k$-tuple 
$\sigma(c) = (\lambda(x_0,x_1), \lambda(x_1, x_2), \dots ,
\lambda(x_{k-1}, x_k))$.  We say that $c$ is an {\it increasing
chain} if the $k$-tuple $\sigma(c)$ is increasing; that is,
if $\lambda(x_0, x_1) \leq \lambda(x_1, x_2) \leq \dots \leq 
\lambda(x_{k-1}, x_k)$.  The edge labeling allows us
to order the maximal chains of any interval of $P$ by
ordering the corresponding $k$-tuples lexicographically.
If $\sigma(c_1)$ lexicographically precedes 
$\sigma(c_2)$ then we say that $c_1$ lexicographically precedes
$c_2$ and we denote this by $c_1  <_L c_2$.

\begin{definition}
An edge labeling is called an {\it EL-labeling} ({\it edge lexicographical
labeling}) if for every interval $[x,y]$ in $P$,
\begin{enumerate}
\item there is a unique increasing maximal chain $c$ in $[x,y]$, and 
\item $c <_L c'$ for all other maximal chains $c^{\prime}$ in $[x,y]$.
\end{enumerate}
\end{definition}

If a graded poset $P$ admits an EL-labeling then its order complex 
is shellable \cite{B1}.
Therefore a 
graded poset that admits an EL-labeling is said to be {\it EL-shellable}.

We now turn our attention to posets that come from CW complexes.
First we review the definition of a {\it regular} CW complex.

\begin{definition}
A CW complex is {\it regular} if the closure
$\overline{C}$ of each cell $C$ is homeomorphic to a closed ball and if 
additionally $\overline{C} \setminus C$ is homeomorphic to a sphere.
\end{definition}

Given a CW complex $\K$, we define its {\it face poset}
$\F(\K)$ to be the set of closed cells ordered by containment
and augmented by a least element $\hat{0}$.  In general, the 
order complex $\Vert \F(\K) - \{\hat{0}\} \Vert$ does not reveal
the topology of $\K$.  However, the following result shows that 
regular CW complexes are combinatorial objects in the
sense that the incidence relation of cells determines their topology.

\begin{proposition} \cite[Proposition~4.7.8]{RedBook}\label{RedTheorem}
Let $\K$ be a regular CW complex.  Then
$\K$ is homeomorphic to 
$\Vert \F(\K) - \{\hat{0}\} \Vert$.
\end{proposition}

It is natural to ask when a poset is the face poset of 
a regular CW complex.  Bjorner \cite{B2} gave a complete answer to this
question, which we now explain.

\begin{deftheorem} \cite{B2} \label{CW}
A poset $P$ is said to be a {\it CW poset} if 
\begin{enumerate}
\item $P$ has a least element $\hat{0}$,
\item $P$ is nontrivial, i.e. has more than one element, and
\item for all $x \in P - \{ \hat{0} \}$, the open interval
  $(\hat{0}, x)$ is homeomorphic to a sphere.
\end{enumerate}
Furthermore, $P$ is a CW poset if and only if it is isomorphic
to the face poset of a regular CW complex.
\end{deftheorem}

Bjorner gave some additional criteria for determining when 
a poset $P$ is a CW poset.  We list below the criterion which will
be most useful to us later.

\begin{proposition}\cite{B2}\label{BjornerCW}
Let $P$ be a nontrivial, finite, pure poset of length $d+1$
with least element $\hat{0}$.  
Let $\hat{P}$ denote the poset $P \cup \{\hat{1}\}$, where $\hat{1}$
is a new greatest element.   
If $\hat{P}$ is shellable and thin then 
$P$ is isomorphic to $\F(\K)$, where $\K$ is a regular 
$d$-dimensional CW complex homeomorphic to the $d$-sphere.
\end{proposition}

A finite graded poset with $\hat{0}$ and $\hat{1}$
is called {\it Eulerian} if every interval of length at least one
has the same number of elements of odd rank as of even rank.
It is an elementary topological result \cite{Stanley} that if
a poset is a CW poset, then it is Eulerian.

\begin{lemma} \cite{Stanley} \label{Eulerian}
A CW poset is Eulerian.
\end{lemma}

\section{Background on the Bruhat order}\label{Bruhat}

In this section we will review properties of the Bruhat order $\leq$
of a Coxeter group, define reflection orders, and 
prove some lemmas that we will need later.
We will assume knowledge of the basic definitions
of Coxeter systems and Bruhat order; we refer the reader to 
\cite{Humphreys}
for details.  

Fix a Coxeter group $W$ generated by a set of simple reflections
$\{ s_i \ \vert \ i \in I \}$ and let $T$ be the set of all reflections. Let
$\ell: W \to \N$ denote the length function: $\ell(w)$ is the length
of a minimal reduced expression for $w$.  Recall that 
length is the rank function for
the Bruhat order of a Coxeter group (which is a graded poset).
One of the key facts about 
the nature of reduced expressions in $W$ is the {\it Strong Exchange
Property}, which is due to Verma \cite{Verma}.

\begin{theorem} (Strong Exchange Property)
Let $w=s_1 \dots s_r$ (where $i \in I$ for $1 \leq i \leq r$), 
not necessarily a reduced 
expression.  Suppose a reflection $t \in T$ satisfies
$\ell(wt) < \ell(w)$.  Then there is an index $i$ for which 
$wt = s_1 \dots \hat{s_i} \dots s_r$ (omitting $s_i$).  If the 
expression for $w$ is reduced, then $i$ is unique.
\end{theorem}

Without loss of generality, assume that $W$ is realized
geometrically as a group of isometries of a real vector space $V$.
Let $\Pi$ denote the set of simple roots, $\Phi^+$ the set of 
positive roots, and for non-isotropic $\alpha \in V$, let 
$r_{\alpha}: V \to V$ denote the reflection in $\alpha$.
Note that the map $\alpha \mapsto r_{\alpha}$ is a bijection 
between $\Phi^+$ and $T$.


Dyer \cite{Dyer} proved that intervals in the Bruhat order of a Coxeter
group are EL-shellable, 
strengthening the earlier CL-shellability result of Bjorner and Wachs 
\cite{BW1}. 
Dyer's primary tool was his notion of ``reflection orders," certain
total orderings of $T$.  To describe reflection orders, we use the 
fact that 
positive roots are naturally in bijection with $T$.  The orders
on the positive roots which correspond to reflection orders on $T$
may be characterized as follows: ``the restriction of the order to the
positive roots lying on the plane spanned by any two positive roots
is one of the two possible orders in which a ray from the origin,
undergoing a full rotation in the plane beginning at a negative root,
would sweep through the positive roots on that plane" \cite{Dyer}.
The 
algebraic definition of a reflection order is the following.

\begin{definition} \cite{Dyer}
A total order $\preceq$ on $T$ is called a {\it reflection order} if for any
dihedral reflection subgroup $W'$ of $W$ either 
$r \prec rsr \prec  rsrsr \prec \dots \prec srsrs \prec srs \prec s$ or 
$s \prec srs \prec srsrs \prec \dots \prec rsrsr \prec rsr \prec r$.
\end{definition}

Dyer proved the existence of many different reflection orders.

\begin{proposition} \cite{Dyer}
Let $J, K$ be disjoint subsets of $I$ and let 
$W_J = \langle  s_j \ \vert \ j \in J  \rangle$,
$W_K = \langle  s_k \ \vert \ k \in K  \rangle$
be the corresponding parabolic subgroups of 
$W$.  Then there is a reflection order $\preceq$ on $T$ such that 
\begin{enumerate}
\item $t \prec t'$ if $t \in W_J \cap T$ and $t' \in T \setminus W_J$, and 
\item $t \prec t'$ if $t \in T \setminus W_K$ and $t' \in W_K \cap T$.
\end{enumerate}
\end{proposition}

In other words, there is a reflection order such that 
reflections contained in
$W_J$ come first, and reflections contained in $W_K$ come last.

Dyer used these reflection orders to prove the following lemma.

\begin{lemma} \cite{Dyer} \label{diamond}
Let $u,w \in W$ with $u \leq w$ and $\ell(w)-\ell(u) = 2$.
Let $\preceq$ be a reflection order. 
Then there exist unique $x, y \in W$ such that 
$u \lessdot x \lessdot w$, $u \lessdot y \lessdot w$, 
$x^{-1}u \prec w^{-1}x$ and
$y^{-1}u \succ w^{-1}y$.  Moreover,
$w^{-1}y \prec w^{-1}x$ and $x^{-1}u \prec y^{-1} u$.
\end{lemma}

In particular, this lemma implies the following result
(originally proved in \cite{BW1}):

\begin{corollary} \cite{BW1, Dyer} The Bruhat order of $W$ is thin.
\end{corollary}

Dyer then extended this result to EL-shellability, as follows.

\begin{proposition} \cite{Dyer} \label{ELDyer}
Fix a reflection order $\preceq$ and 
regard the set $T$ of reflections as a poset under $\preceq$.  
Label
each edge $x \gtrdot y$ of the Bruhat order by the reflection
$x^{-1}y$.  
Then this edge labeling is an EL-labeling; 
therefore the Bruhat order is EL-shellable.
\end{proposition}


We will now prove some further properties of Bruhat order which will be
useful to us later.  Let us fix $(W,I)$ and a subset $J \subset I$.
We will denote by $W^J$ the set of minimal-length coset representatives
of $W/W_J$.  Recall that each element $w \in W$ has a unique factorization
of the form $w' u$, where $w' \in W^J$ and $u \in W_J$.  Also recall
that the projection from $W$ to $W^J$ is
order-preserving.  
Choose a reflection order $\preceq$ such that reflections
contained in $W_J$ come last.
As above, we will label the edge between two elements 
$v \gtrdot u$ with the reflection $v^{-1}u$.





\begin{lemma}\label{useful}
Suppose that $w, w' \in W^J$, and $w'v \lessdot w$ for some $v \in W_J$.
Let $\lambda$ be the label of the edge between $w$ and $w'v$.
Then if $m$ is any chain in $W$ from $w$ to $w'$ which passes through
$w'v$, the labels of $m$ following $\lambda$  are  greater
than $\lambda$.
\end{lemma}

\begin{proof}
Clearly any chain in $W$ from $w$ to $w'$ through $w'v$ will have
the form 
$w \gtrdot w'v \gtrdot w'v_1 \gtrdot w'v_2 \gtrdot \dots 
\gtrdot w'v_{m-1} \gtrdot w'v_m = w'$
where $v_i \in W_J$ and $v_m = e$.  This is because the projection
from $W$ to $W^J$ is order-preserving.  

By the definition of our labeling, $\lambda = w^{-1} w'v$, which is clearly
not in $W_J$.
Now note that the edge between any $w'v_i$ and $w'v_{i+1}$ is labeled by 
the element $(w'v_i)^{-1} w'v_{i+1} = v_i^{-1} v_{i+1}$, which is in $W_J$.
Because we chose $\preceq$ to be a reflection order in which elements 
of $W_J$ come last, the lemma follows.
\end{proof}

This lemma implies the following results.

\begin{corollary}\label{cor1}
Suppose that $w, w' \in W^J$, and $w \gtrdot w'v$ for some $v \in W_J$.
Let $\lambda$ be the label of the edge between $w$ and $w'v$.
Then the unique chain from $w$ to $w'$ in $W$ with an increasing
label 
begins with $\lambda$.
\end{corollary}

\begin{proof}
By Lemma \ref{useful}, 
we know that for any chain from $w$ to $w'$ which goes through 
$w'v$, the edge labels after the initial $\lambda$ are greater than $\lambda$.
Furthermore, by Proposition \ref{ELDyer}, 
there exists a unique increasing chain from $w'v$
to $w'$.  Therefore this gives rise to an increasing chain from 
$w$ to $w'$ whose initial edge label is $\lambda$.  This increasing chain
must necessarily be the unique one.
\end{proof}

\begin{corollary}\label{cor2}
Fix $w, w' \in W^J$ with $w > w'$.  Then there can be at most one 
element in $W$ covered by $w$ which has the form $w'v$
(for any $v \in W_J$).
\end{corollary}

\begin{proof}
Suppose that there were another element covered by $w$ of the
form 
$w'u$ for some $u \in W_J$. 
Let the label of the edge from $w$ to $w'v$ be $\lambda:= w^{-1}w'v$ and 
let the label of the edge from $w$ to $w'u$ be $\mu:= w^{-1}w'u$.
Then Corollary \ref{cor1} implies that the unique increasing chain
in $W$ from $w$ to $w'$ begins with $\lambda$ and also begins with $\mu$, 
which implies that $\lambda=\mu$.  But this implies that $w'v = w'u$.
\end{proof}

\begin{lemma}\label{lem1}
Fix $w',w \in W^J$ with $w' < w$.  Suppose that 
$w'u \lessdot wv$ where $u,v \in W_J$.  Then we can write
$u$ in the form $r v$ where $r \in W_J$ and 
$\ell(u)=\ell(r) + \ell(v)$.
\end{lemma}

\begin{proof}
Let us choose a reduced expression for $wv$, 
say $s_1 s_2 \dots s_m s_{m+1} \dots s_n$ where
$w = s_1 s_2 \dots s_m$ and $v = s_{m+1} \dots s_n$.  Here, the 
$s_i$ are simple reflections.  
Then by the Strong Exchange property, 
$w'u = s_1 s_2 \dots \hat{s}_a \dots s_m s_{m+1} \dots s_n$
for a uniquely determined $a \leq m$.
Clearly $s_1 s_2 \dots \hat{s}_a \dots s_m$ is reduced, 
and hence 
$s_1 s_2 \dots \hat{s}_a \dots s_m \lessdot w$.
We can now write $s_1 s_2 \dots \hat{s}_a \dots s_m$ uniquely
in the form $tr$ where $t \in W^J$ and $r \in W_J$
and $\ell(t)+\ell(r) = m-1$.
But now $w'u = trv$ which implies that $t = w'$ and $u = rv.$
Since $\ell(t) + \ell(r) = m-1$, $\ell(v) = n-m$, and 
$\ell(w'u) = n-1$, we must have that 
$\ell(rv) = \ell(r) + \ell(v)$.  
\end{proof}

\begin{corollary}\label{cor3}
Fix $w',w \in W^J$ with $w' < w$ and fix $v \in W_J$.  Then
$wv$ covers at most one element of the form $w'u$
(for any $u \in W_J$).
\end{corollary}

\begin{proof}
This follows from Corollary \ref{cor2} and Lemma \ref{lem1}.
\end{proof}

\begin{corollary}\label{cor4}
Fix $w', w\in W^J$ with $w' < w$ and fix $u \in W_J$.  Then
$w'u$ is covered by at most one element of the form
$wv$ (for any $v \in W_J$).
\end{corollary}

\begin{proof}
Assume that $w'u$ is covered by two elements $wv$ and $w\tilde{v}$.
Then by Lemma \ref{lem1}, we can write $u = rv$ and also
$u=\tilde{r}\tilde{v}$, where $\ell(u) = \ell(r) + \ell(v)$ and 
$\ell(u) = \ell(\tilde{r}) + \ell(\tilde{v})$.   Since
$wv$ and $w \tilde{v}$ cover $w'rv = w'\tilde{r} \tilde{v}$, it follows
that $w$ covers $w'r$ and $w'\tilde{r}$.  But now by Corollary \ref{cor2},
we must have that $r = \tilde{r}$.  It follows that $v=\tilde{v}$.
\end{proof}



\section{The cell decomposition of $\PP^J_{\geq 0}$ and the 
poset of cells $\Q^J$}\label{Poset}

In this section we introduce the totally nonnegative part of a flag
variety $\PP^J_{\geq 0}$ and its cell decomposition proved by
Rietsch. 
Note that the reader 
interested purely in the combinatorial results of this paper may skip
most of the content in this section, focusing on 
Definition \ref{MainDef}
for the definition of the 
poset $\Q^J$.

Let $G$ be a semisimple linear algebraic group over $\C$ split over
$\R$, with split torus $T$.  We identify $G$ (and related spaces)
with their real points and consider them with their real topology.
Let $X(T) = \Hom(T, \R^*)$ and $\Phi \subset X(T)$ the set of roots.
Choose a system of positive roots $\Phi^+$.  We denote by $B^+$ the
Borel subgroup corresponding to $\Phi^+$ and by $U^+$ its unipotent
radical.  We also have the opposite Borel subgroup $B^-$ such that 
$B^+ \cap B^- = T$, and its unipotent radical $U^-$.

Denote the set of simple roots by 
$\Pi = \{\alpha_i \ \vline \ i \in I \} \subset R^+$.
For each $\alpha_i \in \Pi$ there is an associated homomorphism
$\phi_i : \Sl_2 \to G$.
Consider the $1$-parameter subgroups in $G$ (landing in $U^+, U^-$,
and $T$, respectively) defined by
\begin{equation*}
x_i(m) = \phi_i \left(
                   \begin{array}{cc}
                     1 & m \\ 0 & 1\\
                   \end{array} \right) ,\    
y_i(m) = \phi_i \left(
                   \begin{array}{cc}
                     1 & 0 \\ m & 1\\
                   \end{array} \right) ,\    
\alpha_i^{\vee}(t) = \phi_i \left(
                   \begin{array}{cc}
                     t & 0 \\ 0 & t^{-1}\\
                   \end{array} \right) ,    
\end{equation*}
where $m \in \R, t \in \R^*, i \in I$.
The datum $(T, B^+, B^-, x_i, y_i; i \in I)$ for $G$ is 
called a {\it pinning}.  The standard pinning for 
$\Sl_d$ consists of the diagonal, upper-triangular, and lower-triangular
matrices, along with the simple root subgroups 
$x_i (m) = I_d + mE_{i,i+1}$ and 
$y_i (m) = I_d + mE_{i+1,i}$ where 
$I_d$ is the identity matrix and $E_{i,j}$ has a $1$ in  
position $(i,j)$ and zeroes elsewhere.

The Weyl group $W = N_G(T) / T$ acts on $X(T)$ permuting the roots 
$\Phi$.  The simple reflections $s_i \in W$ are given explicitly
by $s_i:= \dot{s_i} T$ where
$\dot{s_i} := 
                 \phi_i \left(
                   \begin{array}{cc}
                     0 & -1 \\ 1 & 0\\
                   \end{array} \right)$
and any $w \in W$ can be expressed as a product
$w = s_{i_1} s_{i_2} \dots s_{i_m}$ with $\ell(w)$ factors.  We set 
$\dot{w} = \dot{s_{i_1}} \dot{s_{i_2}} \dots \dot{s_{i_m}}$.

We can identify the flag variety with the variety $\B$
of Borel subgroups, via 
\begin{equation*}
gB^+ \Longleftrightarrow g \cdot B^+ := gB^+ g^{-1}.
\end{equation*}
We have the Bruhat decompositions
\begin{equation*}
\B = \sqcup_{w \in W} B^+ \dot{w} \cdot B^+ = 
    \sqcup_{w \in W} B^- \dot{w} \cdot B^-
\end{equation*}
of $\B$ into $B^+$-orbits called {\it Bruhat cells},
and $B^-$-orbits called {\it opposite Bruhat cells}.

\begin{definition}
For $v, w \in W$ define 
\begin{equation*}
R_{v,w}: = B^+ \dot{w} \cdot B^+ \cap B^- \dot{v} \cdot B^+.
\end{equation*}
\end{definition}
The intersection $R_{v,w}$ is non-empty precisely if $v \leq w$,
and in that case is irreducible of dimension $\ell(w) - \ell(v)$,
see \cite{KL}.

Let $J \subset I$.  The parabolic subgroup 
$W_J \subset W$ 
corresponds to a parabolic subgroup $P_J$ in $G$
containing $B^+$.  Namely, 
$P_J = \sqcup_{w \in W_J} B^+ \dot{w} B^+$.
Let $\PP^J$ denote the variety of all parabolic subgroups of $G$
conjugate to $P_J$.  This can be identified with the partial 
flag variety $G/P_J$ via
\begin{equation*}
gP_J \Longleftrightarrow gP_J g^{-1}.
\end{equation*}
We have the usual projection from the full flag variety to 
a partial flag variety which takes the form 
$\pi = \pi^J: \B \to \PP^J$, where 
$\pi(B)$ is the unique parabolic subgroup of type $J$ containing
$B$.

We now give the relevant definitions from total positivity.   
\begin{definition} \cite{Lusztig3}
The totally nonnegative 
part $U_{\geq 0}^-$ of $U^-$ is defined to be the semigroup in 
$U^-$ generated by the $y_i(t)$ for $t \in \R_{\geq 0}$.  

The 
totally nonnegative part of $\B$ is defined by 
\begin{equation*}
\B_{\geq 0} := \overline{ \{u \cdot B^+ \ \vline \ u \in U_{\geq 0}^- \} },
\end{equation*}
where the closure is taken inside $\B$ in its real topology.

The totally nonnegative
part $\PP_{\geq 0}^J$ of a partial flag variety  
$\PP^J$ is defined to be 
$\pi^J (\B_{\geq 0})$.  
\end{definition}

There are natural decompositions of $\B_{\geq 0}$ and 
$\PP_{\geq 0}^J$ which were introduced by Lusztig \cite{Lusztig3, 
Lusztig2}.

\begin{definition}  \cite{Lusztig3}
For $v, w \in W$ with $v \leq w$, let 
\begin{equation*}
R_{v,w ; >0} := R_{v,w} \cap \B_{\geq 0}.
\end{equation*}
\end{definition}

We write $W^J$ (respectively $W^J_{max}$) for the set of minimal
(respectively maximal) length coset representatives of $W/W_J$.

\begin{definition} \cite{Lusztig2} \label{index}
Let $\I^J \subset W^J_{max} \times W_J \times W^J$ be the set of 
triples $(x,u,w)$ with the property that $x \leq wu$.  
Given $(x,u,w) \in \I^J$, 
we define $P_{x,u,w; >0}^J := \pi^J(R_{x,wu; >0}) = \pi^J(R_{xu^{-1},w; >0}).$
\end{definition}

Rietsch \cite{Rietsch1} proved that 
$R_{v,w; >0}$ and 
$P_{x,u,w; >0}^J$ are semi-algebraic cells of dimension
$\ell(w)-\ell(v)$ and $\ell(wu) - \ell(x)$, respectively.
Furthermore, 
Rietsch \cite{Rietsch2} has given a concrete description of the order
relation on cells.

\begin{theorem}\cite{Rietsch2} \label{RTheorem}
We have that $P_{x,u,w;>0}^J \subset \overline{P_{x',u',w';>0}^J}$ precisely
if there exist $u_1, u_2 \in W_J$ with $u_1 u_2 = u$ and 
$\ell(u) = \ell(u_1)+\ell(u_2)$, such that the following holds:
\begin{equation*}
x'{u'}^{-1} \leq x u_2^{-1} \leq w u_1 \leq w'.
\end{equation*}
Furthermore, the closure of a cell is a union of cells.
\end{theorem}

\begin{definition}\label{MainDef}
Fix a Coxeter system $(W,I)$ and a subset $J \subset I$ such that 
$W_J$ is finite.
Let $\Q^J$ denote the poset with elements
$\{\hat{0}\} \cup 
\{Q_{x,u,w} \ \vert \ (x,u,w) \in \I^J \}$ and 
with the order relation above.  
That is, 
$Q_{x,u,w} < Q_{x',u',w'}$ if and only 
if there exist $u_1, u_2 \in W_J$ with $u_1 u_2 = u$ and 
$\ell(u) = \ell(u_1)+\ell(u_2)$, such that the following holds:
\begin{equation*}
x'{u'}^{-1} \leq x u_2^{-1} \leq w u_1 \leq w'.
\end{equation*}
Additionally, we stipulate that $\hat{0} < Q_{x,u,w}$ for all $Q_{x,u,w} \in \Q^J$.
\end{definition}

Clearly when $W$ is a Weyl group of a semisimple linear algebraic group,
$\Q^J$ is the poset of closed cells of $\PP^J_{\geq 0}$ augmented by a least
element.  However, we define $\Q^J$ as above because our combinatorial 
results are true in this generality; we thank Anders Bjorner for pointing
this out to us.

Besides having a unique least element $\hat{0}$, $\Q^J$ also has a unique greatest
element:
this element is
$Q_{u_0, u_0, w_0^J}$, where $u_0$ is the longest element in $W_J$ and 
$w_0^J$ is the longest element in $W^J$.  
In Section \ref{EL} we will prove that $\Q^J$ is a graded poset,
where the rank of $Q_{x,u,w}$ is simply $\ell(wu)-\ell(x)$.


See Example \ref{G24} and Figure \ref{Picture} for the case
of the Grassmannian $Gr_{2,4}(\R)$.  

\begin{remark}
When $\PP^J$ is the (type A) Grassmannian, the poset $\Q^J$ 
is the poset of cells of the totally nonnegative Grassmannian, 
studied by Postnikov \cite{Postnikov} and the author \cite{Williams}.
This poset is also called the {\it cyclic Bruhat order},
and Postnikov \cite{Postnikov} has shown that it can be described in terms
of many different combinatorial objects, such as decorated permutations
and $\Le$-diagrams (certain $0-1$ tableaux).  See Appendix \ref{A} for 
more details.
\end{remark}

\begin{remark}
When $\PP^J$ is the complete flag variety, the poset $\Q^J$ has an 
especially simple description.  Recall that if $P$ is a poset, then
the {\it interval poset} $Int(P)$ is defined to be the poset of 
intervals $[x,y]$ of $P$, ordered by containment.  In this 
special case, $\Q^J$ is simply the interval poset of the Bruhat order.
\end{remark}

\section{$\Q^J$ is thin}\label{thin}

In this section we will study cover relations in $\Q^J$ as well
as intervals of rank $2$.  In particular, we will prove that 
$\Q^J$ is thin.

Recall that $Q_{x',v,w'} < Q_{x,u,w}$ if and only if there exist
$v_1$ and $v_2$ such that $v=v_1 v_2$ (lengths add) and 
$x u^{-1} \leq x' v_2^{-1} \leq w'v_1 \leq w$.  Here 
$w,w'\in W^J, u,v, v_1, v_2\in W_J$, and $x,x'\in W^J_{max}$.  There are three
types of cover relations:

\begin{description}
\item[Type 1] $Q_{x',v,w'} \lessdot Q_{x,u,w}$ such that
     $x u^{-1}=x' v_2^{-1}$ and $w'v_1 \lessdot w$.  
  This implies that $x=x', u=v_2$, $w'v \lessdot wu$, 
  and hence $w' < w$.
\item[Type 2] $Q_{x',v,w'} \lessdot Q_{x,u,w}$ such that
     $x u^{-1} \lessdot x' v_2^{-1}$ and $w'v_1 = w$.
  This implies that $w'=w$ and $xu^{-1} \lessdot x' v^{-1}$.  
\item[Type 3]  $\hat{0} \lessdot Q_{x,u,w}$ where $Q_{x,u,w}$
    is a $0$-cell.  This implies that $x=wu$.
\end{description}

We will now prove that intervals of rank $2$ in $\Q^J$ are diamonds.

\def\thetheoremB{\ref{Theorem1}}
\begin{theoremB}
$\Q^J$ is thin.
\end{theoremB}

\begin{proof}
Let us fix two elements $Q_{a,b,c} < Q_{x,u,w} \in \Q^J$, such that 
$\rank(Q_{x,u,w}) - \rank(Q_{a,b,c}) = 2$.  We must show that there
are exactly two elements between $Q_{x,u,w}$ and $Q_{a,b,c}$ in 
$\Q^J$.

Since $Q_{a,b,c} < Q_{x,u,w}$, there exist $b_1, b_2\in W_J$ with 
$b_1 b_2 = b$ (and the lengths add) such that 
\begin{equation} \label{order}
x u^{-1} \leq ab_2^{-1} \leq cb_1 \leq w.
\end{equation}
Since the difference in rank between $Q_{a,b,c}$ and $Q_{x,u,w}$ is 
$2$, there are three possibilities for the inequalities above:
\begin{description}
\item[Case 1] $xu^{-1} = ab_2^{-1}$, and $cb_1 < w$ 
where $\ell(w)-\ell(cb_1) = 2$. 
\item[Case 2] $xu^{-1} \lessdot ab_2^{-1}$ and $cb_1\lessdot w$. 
\item[Case 3] $xu^{-1} < ab_2^{-1}$ where $\ell(ab_2^{-1})-\ell(xu^{-1})=2$, and $cb_1= w$. 
\end{description}

Let us first consider Case 1. 
Since $xu^{-1} = ab_2^{-1}$, we have that $x=a$ and $u = b_2$.
The inequality $cb_1<w$ 
implies that $c < w$.
Finally, $cb_1 < w$ and $b_2 = u$ implies that 
$cb < wu$ where $\ell(wu) - \ell(cb)=2$.  

The fact that the Bruhat order is thin implies that there exist
two elements in $W$ between $cb$ and $wu$.  Let us factor these
two elements (uniquely) as $tr$ and $t'r'$ where $t,t'\in W^J$
and $r,r'\in W_J$.  

We now claim that 
\begin{equation}\label{chain1}
Q_{a,b,c} \lessdot Q_{a,r,t}.
\end{equation} and 
\begin{equation}\label{chain2}
Q_{a,r,t} \lessdot Q_{a,u,w}.
\end{equation} 

(The same relations will hold for $Q_{a,r',t'}$.)

First note that $Q_{a,r,t} \in \Q^J$, since 
$a \leq cb < tr$ implies that $a \leq tr$.  

Let us prove (\ref{chain1}).
Note that either $c = t$ or $c < t$.  
If $c=t$, then $cb \lessdot tr$ implies that $b \lessdot r$.
Clearly $a r^{-1} \lessdot ab^{-1} \leq c = c$ which implies
that $Q_{a,b,c} \lessdot Q_{a,r,c}$.  

If $c < t$, then 
since $cb \lessdot tr$, Lemma \ref{lem1} implies that 
we can write $b = b_1 r$ where $\ell(b) = \ell(b_1)+\ell(r)$.  
Let $b_2 = r$.   Since $cb \lessdot tr = tb_2$, it follows that
$cb_1 \lessdot t$.  We now have 
$a r^{-1} = a b_2^{-1} \leq cb_1 \lessdot t$, 
which implies that $Q_{a,b,c} \lessdot Q_{a,r,t}$.

Now let us prove (\ref{chain2}).  
Note that either $t = w$ or $t < w$.  
If $t=w$ then 
$tr \lessdot wu$ implies that 
$r \lessdot u$.   
Therefore $a u^{-1} \lessdot ar^{-1} \leq w = w$ which implies that 
$Q_{a,r,t}=Q_{a,r,w} \lessdot Q_{a,u,w}$.

If $t < w$, then since $tr \lessdot wu$, Lemma \ref{lem1} implies that
we can write $r = r_1 u$ where $\ell(r) = \ell(r_1)+ \ell(u)$.  
Let $r_2 = u$.
Now $tr \lessdot wu$ implies that $tr_1 u \lessdot wu$ which implies
that $tr_1 \lessdot w$.  
We now have $au^{-1} = a r_2^{-1} \leq tr_1 \lessdot w$, and so
$Q_{a,r,t} \lessdot Q_{a,u,w}$.

We have now shown that there exist two elements $Q_{a,r,t}$ and 
$Q_{a,r',t'}$ which lie between $Q_{a,b,c}$ and $Q_{x,u,w}=Q_{a,u,w}$
when we are in the situation of Case 1.  
To complete the proof in this case, we 
need to show that these are the {\it only} two elements which lie
in this open interval.

Take any element $Q_{d,h,f}$ such that 
$Q_{a,b,c} \lessdot Q_{d,h,f} \lessdot Q_{a,u,w}$.  
By using the fact that the projection from $W$ to $W^J$ is order preserving,
it is easy to see that we must have $d=a$.  Now the fact that 
$Q_{a,b,c} \lessdot Q_{a,h,f} \lessdot Q_{a,u,w}$ implies that 
there exist $b_1, b_2\in W_J$ with $b_1 b_2 = b$ such that 
\begin{equation}\label{ineq1}
a h^{-1} \leq ab_2^{-1} \leq cb_1 \leq f,
\end{equation} 
and there exist
$h_1, h_2\in W_J$ with $h_1 h_2 = h$ such that 
\begin{equation}\label{ineq2}
a u^{-1} \leq a h_2^{-1} \leq f h_1 \leq w.  
\end{equation}
The fact that (\ref{ineq1}) and (\ref{ineq2}) represent {\it cover}
relations in the poset implies that in each of (\ref{ineq1}) and (\ref{ineq2})
exactly one of the outer $\leq$'s is a $\lessdot$ and one is an equality.
If $a h^{-1} = ab_2^{-1}$ and $cb_1 \lessdot f$ then
$h = b_2$ and $cb_1 \lessdot f$ implies that $cb \lessdot fh$.  
On the other hand, if 
$a h^{-1} \lessdot ab_2^{-1}$ and $cb_1 = f$ then
$c=f$ and $b_1 = e$.  Thus
$a h^{-1} \lessdot ab^{-1}$ which implies that $b \lessdot h$
and hence $cb \lessdot fh$.
Similarly, it is easy to show from (\ref{ineq2}) that 
$fh \lessdot wu$.  
But now since $cb \lessdot fh \lessdot wu$, 
the element $Q_{a,h,f}$ must be one of the two 
elements $Q_{a,r,t}$, $Q_{a,r',t'}$ which we already found.

Now we consider Case 2.
Recall that we have 
$x u^{-1} \lessdot ab_2^{-1}$ and $cb_1 \lessdot w$, where
$b = b_1 b_2$ and the lengths add.
If $x=a$, then as in the previous case, we will have 
$cb < wu$ with $\ell(wu)-\ell(cb) = 2$, and the two 
elements between $Q_{a,b,c}$ and $Q_{a,u,w}$ will be of the form
$Q_{a,r,t}$ where $cb \lessdot tr \lessdot wu$.  

We now assume that $x\neq a$.
We claim that the only elements between $Q_{a,b,c}$ and 
$Q_{x,u,w}$ are $Q_{a,b_2,w}$ and $Q_{x,b_1 u, c}$.  

To show that $Q_{a,b_2,w}$ is well-defined (i.e.\ that $(a,b_2,w)$ is in the 
indexing set for Rietsch's cell decomposition), note that 
$cb_1 \lessdot w$ implies that $cb \lessdot wb_2$.  This together
with $a \leq cb$ implies that $a \leq wb_2$, and so 
$Q_{a,b_2,w}$ is well-defined.  Now it is easy to see that 
$ab_2^{-1} = ab_2^{-1} \leq cb_1 \lessdot w$ which implies that 
$Q_{a,b,c} \lessdot Q_{a,b_2,w}$, and also 
$xu^{-1} \lessdot ab_2^{-1} \leq w = w$ which implies that 
$Q_{a,b_2,w} \lessdot Q_{x,u,w}$.

In order to show that $Q_{x,b_1 u,c}$ is well-defined, we must first
observe that $\ell(b_1 u) = \ell(b_1) + \ell(u)$. 
By applying Lemma \ref{lem1} to the fact that $xu^{-1} \lessdot ab_2^{-1}$
(and remembering that each element of $W^J_{max}$ has a unique expression
as a product of an element of $W^J$ and the longest element $u_0 \in W_J$),
we see that we can write $b_2$ in the form $u b'$ where 
$\ell(b_2) = \ell(u)+\ell(b')$.  
But now since lengths add in the product $b = b_1 b_2$, this implies 
that lengths add in the product $b_1 u$.  
Therefore $x \leq cu$ implies that $x \leq c b_1 u$, and hence
$Q_{x,b_1 u, c}$ is well-defined.

Now $x u^{-1} \lessdot ab_2^{-1}$ together with the fact that lengths
add in the product $b_1 u$ implies that 
$x u^{-1} b_1^{-1} \lessdot ab_2^{-1}b_1^{-1}$.  Therefore
$x (b_1 u)^{-1} \lessdot a b^{-1} \leq c = c$ and so 
$Q_{a,b,c} \lessdot Q_{x,b_1 u, c}$.  
Finally, we have $x u^{-1} = x u^{-1} \leq cb_1 \lessdot w$ 
and so $Q_{x, b_1 u, c} \lessdot Q_{x,u,w}$.  

To complete the proof for Case 2, we need to show that 
the {\it only} two elements between $Q_{a,b,c}$ and 
$Q_{x,u,w}$ are $Q_{a,b_2,w}$ and $Q_{x,b_1 u, c}$.  

Consider any $Q_{r,s,t}$ such that 
$Q_{a,b,c} \lessdot Q_{r,s,t} \lessdot Q_{x,u,w}$.
This implies that there exist $b_1, b_2 \in W_J$ with $b_1 b_2 = b$
and $s_1, s_2 \in W_J$ with $s_1 s_2 = s$ such that 
$rs^{-1} \leq a b_2^{-1} \leq cb_1 \leq t$ and 
$xu^{-1} \leq rs_2^{-1} \leq ts_1 \leq w$.
Since these inequalities represent cover relations in the poset, 
we see that if 
$t = c < w$ then $ts_1 \lessdot w$ and hence $xu^{-1} = rs_2^{-1}$ which
implies that $r = x$.
On the other hand, if 
$t > c$, then we must have $r s^{-1} = ab_2^{-1}$ which implies that $r=a$.  
Note that if $t > c$ then we must in fact have $t=w$, because
as we saw before, $t < w$ 
implies that $r = x$.  Since $r=a$ and we have assumed that $x \neq a$,
this is a contradiction.

Therefore we know that elements between $Q_{a,b,c}$ and $Q_{x,u,w}$
either have the form $Q_{x,s,c}$ or $Q_{a,s,w}$.  
If $Q_{a,b,c} \lessdot Q_{x,s,c} \lessdot Q_{x,u,w}$ then
$xu^{-1} \leq xs_2^{-1} \leq cs_1 \leq w$ together with the fact that 
$c \neq w$ implies that 
$cs_1 \lessdot w$ and $xu^{-1} = xs_2^{-1}$.  
Corollary \ref{cor2} applied to $cs_1 \lessdot w$ 
implies that $s_1$ is uniquely 
determined, and since $s_2 = u$, we have that $s$ is uniquely 
determined.
Therefore there is at most one element of the form $Q_{x,s,c}$
between $Q_{a,b,c}$ and $Q_{x,u,w}$

Similarly, if 
$Q_{a,b,c} \lessdot Q_{a,s,w} \lessdot Q_{x,u,w}$ then
$xu^{-1} \leq as_2^{-1} \leq ws_1 \leq w$ implies that 
$xu^{-1} \lessdot as_2^{-1}$ and $s_1 = e$.  
Corollary \ref{cor4} implies that $s_2$ is uniquely
determined, and hence $s$ is uniquely determined.  
Therefore there is at most one element of the form $Q_{a,s,w}$
between $Q_{a,b,c}$ and $Q_{x,u,w}$.
This completes the proof for Case 2.

The proof for Case 3 is the simplest of all.
The fact that $cb_1  = w$ implies that $b_1 = e$ and $c=w$.
Therefore $xu^{-1} < ab^{-1}$ and 
$\ell(ab^{-1})-\ell(xu^{-1})=2$.  Since the Bruhat order is thin,
there are exactly two elements in the open interval
$(xu^{-1}, ab^{-1})$, which we can factor uniquely as 
$rs^{-1}$ and $r's'^{-1}$ for $r, r' \in W_{max}^J$ and 
$s, s' \in W_J$.  It is clear that 
$Q_{a,b,w} < Q_{r,s,w} < Q_{x,u,w}$ because 
$rs^{-1} < ab^{-1} \leq w \leq w$ and 
$xu^{-1} < rs^{-1} \leq w \leq w$.  Conversely, it is easy to see that if 
some $Q_{\tilde{r}, \tilde{s}, w}$ satisfies 
$Q_{a,b,w} < Q_{\tilde{r},\tilde{s},w} < Q_{x,u,w}$, then
$xu^{-1} \leq \tilde{r} \tilde{s}^{-1} \leq ab^{-1}$.

To complete the proof,
we must address the rank $2$ intervals whose least element is
$\hat{0}$.  Let the greatest element of such an interval be
$Q_{x,u,w}$.  It follows that $x \lessdot wu$.  It is now an easy
exercise to see that there are exactly two elements in the open
interval $(\hat{0}, Q_{x,u,w})$: 
$Q_{wu_0, u_0, w}$ and $Q_{x, u_0, xu_0}$, where $u_0$ is the 
longest element in $W_J$.
Therefore $\Q^J$ is thin.
\end{proof}

We now summarize the analysis of the previous proof.

\begin{remark}
In the situation of Case 1, the interval
$[Q_{x,b,c}, Q_{x,u,w}]$ in $\Q^J$ naturally corresponds to the 
interval $[cb, wu]$ in $W$.  
Note that because $c < w$, each chain from $Q_{x,u,w}$ to $Q_{x,b,c}$ must
contain at least one Type 1 cover relation; however, a Type 2 
cover relation may also occur in this chain.

In the situation of Case 2, 
one of the chains from $Q_{x,u,w}$ 
to $Q_{a,b,c}$ consists of a Type 1 cover relation followed by 
a Type 2 cover relation; the other chain consists of a Type 2 cover
relation followed by a Type 1 cover relation.

In the situation of Case 3, the interval
$[Q_{a,b,w}, Q_{x,u,w}]$ in $\Q^J$ naturally corresponds to the 
interval $[xu^{-1}, ab^{-1}]$ in $W$.  Note that all four edges
of the diamond interval 
$[Q_{a,b,w}, Q_{x,u,w}]$ must
have Type 2.

For rank $2$ intervals of the form $[\hat{0}, Q_{x,u,w}]$, 
one of the chains from $Q_{x,u,w}$ to $\hat{0}$ has the form
Type 1 -- Type 3, and the other chain has the form
Type 2 -- Type 3.
\end{remark}

\section{$\Q^J$ is EL-shellable}\label{EL}

In this section we will prove that $\Q^J$ is a graded poset and that 
it has an 
EL-labeling, which implies that $\Q^J$ is EL-shellable \cite{B1}.
In particular,
the order complex of $\Q^J$ is shellable.

Let us begin by  
recalling the three types of cover relations in $\Q^J$:

\begin{description}
\item[Type 1] $Q_{x',v,w'} \lessdot Q_{x,u,w}$ such that
     $x u^{-1}=x' v_2^{-1}$ and $w'v_1 \lessdot w$. 
  This implies that $x=x', u=v_2$, and $w'v \lessdot wu$, and hence
  $w' < w$.
\item[Type 2] $Q_{x',v,w'} \lessdot Q_{x,u,w}$ such that
     $x u^{-1} \lessdot x' v_2^{-1}$ and $w'v_1 = w$.
  This implies that $w'=w$ and $xu^{-1} \lessdot x' v^{-1}$.  
\item[Type 3]  $\hat{0} \lessdot Q_{x,u,w}$ where $Q_{x,u,w}$
    is a $0$-cell.  This implies that $x=wu$.
\end{description}

We now prove a lemma  
which describes a condition that diamond intervals
in $\Q^J$ may {\it not} possess.

\begin{lemma}\label{forbidden}
There are no diamond intervals in $\Q^J$ 
in which the top two edges have Type 2 and
the bottom two edges have Type 1.
\end{lemma}

\begin{proof}
Let $P_1, P_2, P_2^{\prime}, P_3$ be the elements of a diamond interval in $\Q^J$, such that
$P_3 \lessdot P_2 \lessdot P_1$ and also $P_3 \lessdot P_2^{\prime} 
\lessdot P_1$.  
%
The element $P_1$ has the form
$Q_{x,u,w}$ and if the top two edges have Type 2, then
we can write $P_2 = Q_{x',u',w}$ and $P_2^{\prime} = Q_{x'',u'',w}$.
If additionally the bottom two edges have Type 1 then
$P_3 = Q_{\tilde{x},\tilde{u},\tilde{w}}$ where
$\tilde{x} = x'$ and $\tilde{x} = x''$.  Therefore $x'=x''$.
Furthermore, 
$\tilde{w}\tilde{u} \lessdot w u'$ and
$\tilde{w}\tilde{u} \lessdot w u''$, which implies by 
Corollary \ref{cor4} that $u' = u''$.  But this shows that 
$P_2 = P_2^{\prime}$, which contradicts the fact that we 
were considering a diamond interval.
Therefore the kind of diamond interval described 
in Lemma \ref{forbidden}
is impossible.
\end{proof}

\begin{proposition}\label{prop1}
Suppose that $Q_{a,b,c} < Q_{x,u,w}$ where $c<w$.  Then 
there exists some $Q_{r,s,t}$ with $t < w$ such that 
$Q_{a,b,c} < Q_{r,s,t} \lessdot Q_{x,u,w}$.  
\end{proposition}

\begin{proof}
We prove this by induction on $\ell(w) - \ell(c)$.  

The base case is when $c \lessdot w$.  
In this situation we claim that 
$Q_{a,b,c} < Q_{x,u,c} \lessdot Q_{x,u,w}$.  
The fact that $Q_{a,b,c} < Q_{x,u,w}$ implies that there exists a 
decomposition $b = b_1 b_2$ (lengths add) such that 
$x u^{-1} \leq a b_2^{-1} \leq cb_1 \leq w$.  
Here we cannot have $cb_1 = w$ because that would imply that $c=w$.
Therefore $cb_1 < w$ and the fact that $c \lessdot w$ implies that 
$b_1 = e$.  
This implies that $x u^{-1} \leq ab^{-1} \leq c \lessdot w$.
In particular, $xu^{-1} \leq c$ implies that $Q_{x,u,c}$ is well-defined.
Now it is obvious that $Q_{x,u,c} \lessdot Q_{x,u,w}$ since
$xu^{-1} \leq xu^{-1} \leq c \lessdot w$.  
It remains to show that $Q_{a,b,c} < Q_{x,u,c}$, i.e. that 
there exist $\tilde{b}_1, \tilde{b}_2$ such that $b =
\tilde{b}_1\tilde{b}_2$ (lengths add) and  
$x u^{-1} \leq a \tilde{b}_2^{-1} \leq c \tilde{b}_1 \leq c$.
Clearly we can take $\tilde{b}_1 = e$ and $\tilde{b}_2 = b$.

We now prove the general case.  Consider the interval $I$ between
$Q_{a,b,c}$ and $Q_{x,u,w}$.  Note that any $Q_{r,s,t}$ that lies 
in this interval necessarily satisfies $c \leq t \leq w$.  
Now if any $Q_{r,s,t}$ satisfies $c < t < w$ then since
$\ell(w)-\ell(t) < \ell(w) - \ell(c)$, we are done by induction.
Therefore we are left to consider the case that for each $Q_{r,s,t}$
in this interval, either $t=c$ or $t=w$.  Let us choose a 
$Q_{r,s,c}$ in this interval with maximal rank.  If $Q_{r,s,c}$ is 
covered by $Q_{x,u,w}$, then we are done.  If not, then all elements 
$Q_{e,f,g}$ of $I$ which are greater than $Q_{r,s,c}$ satisfy 
$g=w$.  In particular, there is a diamond interval with $Q_{r,s,c}$ 
at the bottom in which the other three elements have the form
$Q_{*,*,w}$.  But this is impossible by Lemma \ref{forbidden}.
\end{proof}

\begin{lemma}\label{prop2}
Suppose that $Q_{a,b,w} < Q_{x,u,w}$.  Then there exists some 
$Q_{r,s,w}$ such that 
$Q_{a,b,w} < Q_{r,s,w} \lessdot Q_{x,u,w}$.  
\end{lemma}

\begin{proof}
We have $xu^{-1} \leq ab_2^{-} \leq wb_1 \leq w$ which implies that 
$xu^{-1} \leq ab^{-1} \leq w$.  Because the Bruhat order is graded, 
there exists an element $v$ such that $xu^{-1} \lessdot v \leq ab^{-1}$
which can be factored uniquely in the form $v = rs^{-1}$ where 
$r \in W^J_{max}, s \in W_J$.  It is now easy to see that 
$Q_{r,s,w} \in \Q^J$ and that $Q_{a,b,c} < Q_{r,s,w} \lessdot Q_{x,u,w}$.
\end{proof}

\begin{corollary}
$\Q^J$ is a graded poset, where the rank of $Q_{x,u,w} = \ell(wu)-\ell(x)$.
\end{corollary}

\begin{proof}
This follows from Proposition \ref{prop1} and Lemma \ref{prop2}.
\end{proof}

We now propose an edge labelling for $\Q^J$.  
See Example \ref{G24} together with 
Figure \ref{Picture} for the example of the
Grassmannian $Gr_{2,4}(\R)$.

\begin{definition}
Label Type 1 edges with the element $(wu)^{-1} w'v \in T$; 
label Type 2 edges with the element $(x'v^{-1})^{-1}xu^{-1} \in T$;
and label Type 3 edges with the symbol $\emptyset$.
Choose any reflection order $\preceq$ on $T$ such that elements of $T\cap W_J$
come last.  We then choose the total order $\LL$ on 
labels of edges in $\Q^J$ determined by the following conditions:
\begin{enumerate}
\item If $\lambda$ is any Type 1 label and $\mu$ is any Type 2 label, then 
  $\lambda \LL \emptyset \LL \mu$.
\item If $\lambda$ and $\mu$ are Type 1 labels then $\lambda \LL \mu$ if and only 
if $\lambda \prec \mu$.
\item If $\lambda$ and $\mu$ are Type 2 labels then $\lambda \LL \mu$ if and only if
$\lambda \prec \mu$.
\end{enumerate}
\end{definition}

\begin{remark} \label{crucial}
Observe that the labels of Type 1 edges are never in $W_J$.  This will 
be important for our arguments later.
\end{remark}

\def\thetheoremB{\ref{Theorem2}}
\begin{theoremB}
The labeling of edges of $\Q^J$ described above is an EL-labeling.
\end{theoremB}

\begin{proof}
Fix two elements $Q_{x,u,w}$ and $Q_{a,b,c}$ in $\Q^J$, such that 
$Q_{a,b,c} < Q_{x,u,w}$.
First we will show that the lexicographically minimal chain (with 
respect to 
$\LL$) from $Q_{x,u,w}$ to 
$Q_{a,b,c}$ is increasing.

By Proposition \ref{prop1}, it is clear that the lexicographically
minimal chain $m$ will consist of a series of Type 1 edges followed by 
Type 2 edges; let the chain label be
$(\lambda_1, \lambda_2, \dots, \lambda_i, 
\mu_1, \mu_2, \dots, \mu_j)$, where the $\lambda_k \in T$ are Type 1 edges
and the $\mu_k \in T$ are Type 2 edges.

For the sake of contradiction, 
suppose that for some $h$ we have $\lambda_h \RR \lambda_{h+1}$.
In other words, 
$\lambda_h \succ \lambda_{h+1}$,
where $\lambda_h$ and $\lambda_{h+1}$ are the labels, 
respectively, for the edges of the chain 
$Q_{x,u_1, w_1} \gtrdot Q_{x,u_2, w_2} \gtrdot Q_{x,u_3,w_3}$.
As shown in Section \ref{thin}, the interval from 
$Q_{x,u_1, w_1}$ to $Q_{x,u_3, w_3}$ is a diamond.  If we let 
$Q_{x,u_2^{\prime}, w_2^{\prime}}$ denote the other middle element of 
this interval, then we know that
$w_3 u_3 \lessdot w_2 u_2 \lessdot w_1 u_1$ and also 
$w_3 u_3 \lessdot w_2^{\prime} u_2^{\prime} \lessdot w_1 u_1$. 
Since $\lambda_h$ and $\lambda_{h+1}$ are Type 1 labels, we 
know that $w_3 < w_2 < w_1$.  However, for $w_2^{\prime}$, 
we know only that $w_3 \leq w_2 \leq w_1$.
Observe that the labels $\lambda_h$ and $\lambda_{h+1}$ that we've used
for our elements in $\Q^J$ are also the labels used by Dyer for 
the edges of the interval $[w_3 u_3, w_1 u_1]$.
Therefore Lemma \ref{diamond} implies that the labels
$\gamma_h$ and $\gamma_{h+1}$ for the 
chain $w_1 u_1 \gtrdot w_2^{\prime} u_2^{\prime} \gtrdot 
w_3 u_3$ satisfy
$\gamma_h \prec \gamma_{h+1}$ and also 
$\gamma_{h} \prec \lambda_h$.
We now claim that the 
chain $(\lambda_1, \lambda_2, \dots , \gamma_h, \gamma_{h+1}, \dots, 
\lambda_i, \mu_1, \mu_2, \dots, \mu_j)$ is lexicographically 
smaller than 
$(\lambda_1, \lambda_2, \dots, \lambda_i, 
\mu_1, \mu_2, \dots, \mu_j)$.  To complete the proof of the claim,
it suffices to show that 
$\gamma_h$ is the label of a {\it Type 1} cover relation.
It is not hard to see that $\gamma_h$ is the label of a Type 1 
cover relation if and only if $\gamma_h \notin W_J$.
Since $w_3 < w_1$, it is clear that we cannot have both 
$\gamma_h$ and $\gamma_{h+1}$ in $W_J$; at least one is {\it not} in 
$W_J$.
And now because we've chosen a reflection order in which elements 
of $W_J$ come {\it last}, it follows that $\gamma_h$ is not in $W_J$,
and hence is the label of a Type 1 cover relation.
We've now found a lexicographically smaller chain, which is a contradiction.

Therefore the lexicographically minimal chain label
$(\lambda_1, \lambda_2, \dots, \lambda_i, 
\mu_1, \mu_2, \dots, \mu_j)$
satisfies $\lambda_1 \LL \lambda_2 \LL \dots \LL \lambda_i$.
Let $Q_{r,s,c}$ denote the element 
of $\Q^J$ that we reach if we start at $Q_{x,u,w}$ and traverse the 
edges labeled by $\lambda_1, \lambda_2, \dots , \lambda_i$.
We need to show that the lexicographically minimal 
chain from $Q_{r,s,c}$ to $Q_{a,b,c}$ is increasing.  First note that all edge 
labels of this interval are Type 2 labels.
Furthermore, by considering the order relation in $\Q^J$, 
it is easy to see that any element $Q_{d,e,c}$ in the interval
$[Q_{a,b,c}, Q_{r,s,c}]$ satisfies $rs^{-1} \leq de^{-1} \leq ab^{-1}$.
Conversely,
for any $d \in W^J_{max}$ and $e \in W_J$ such that 
$rs^{-1} \leq de^{-1} \leq ab^{-1}$, we have that 
$Q_{a,b,c} \leq Q_{d,e,c} \leq Q_{r,s,c}$.  
Therefore the interval $[Q_{a,b,c}, Q_{r,s,c}]$ is isomorphic to the 
dual of the interval $[rs^{-1}, ab^{-1}]$.  Recall that our 
edge labeling 
of $[Q_{a,b,c}, Q_{r,s,c}]$ is inherited from (the dual of)
Dyer's edge labeling of  
$[rs^{-1}, ab^{-1}]$.  It follows -- using EL-shellability of intervals
in Bruhat order -- that the lexicographically minimal
chain (with respect to $\LL$) from $Q_{r,s,c}$ to $Q_{a,b,c}$ is increasing.
Since we've chosen an ordering in which Type 1 labels 
precede Type 2 labels, we have now shown that 
the lexicographically minimal 
chain from $Q_{x,u,w}$ to $Q_{a,b,c}$ is increasing.

It remains to show that this is the {\it unique} increasing chain
from $Q_{x,u,w}$ to $Q_{a,b,c}$.
Clearly the labels on any increasing chain will again consist of a series
of Type 1 edge labels followed by a series of Type 2 edge labels.
Suppose that there are two increasing chains $m_1$ and $m_2$ from $Q_{x,u,w}$
to $Q_{a,b,c}$; let $Q_{r,v,c}$ and $Q_{r', v', c}$ be the 
two intermediate elements of these chains which we obtain after starting
at $Q_{x,u,w}$ and traversing the Type 1 edges of $m_1$ and $m_2$, 
respectively.
It follows that $cv < wu$ and $cv' < wu$.  Furthermore,
the increasing chain labels from $Q_{x,u,w}$ to $Q_{r,v,c}$ and 
from $Q_{x,u,w}$ to $Q_{r',v',c}$ correspond to 
increasing chain labels from $wu$ to $cv$ and from $wu$ to 
$cv'$.  Note that by Remark \ref{crucial}, the labels on these
increasing chains from $wu$ to $cv$ and from $wu$ to $cv'$ are {\it not}
in $W_J$.

Now consider the interval $[c,wu]$.  Clearly both $cv$ and $cv'$ are 
in this interval.  By EL-shellability of the Bruhat order, 
we can find increasing chains from $cv$ to $c$ and also from 
$cv'$ to $c$.  Clearly the labels on these chains will be elements
of $T \cap W_J$.  And therefore by our choice of reflection ordering,
the increasing chain from $wu$ to $cv$ extends to an increasing
chain from $wu$ to $c$; similarly, 
the increasing chain from $wu$ to $cv'$ extends to an increasing
chain from $wu$ to $c$.  We have now found two increasing chains
from $wu$ to $c$, which contradicts the fact that reflection orders
give EL-labelings of the Bruhat order.
Therefore there is a {\it unique} increasing chain from 
$Q_{x,u,w}$ to $Q_{a,b,c}$.  

This would complete the proof that $\Q^J$ is EL-shellable except that
we have so far ignored the chains from $Q_{x,u,w}$ to $\hat{0}$.
We will now address these chains.

Consider all maximal chains from $Q_{x,u,w}$ to a $0$-cell, 
i.e.\ an element $Q_{a,b,c}$ with $\ell(cb)=\ell(a)$.  We claim
that among these, the lexicographically minimal chain $m$ consists 
entirely of Type 1 edges and is increasing.  First note that if 
any element $Q_{a,b,c}$ is not a $0$-cell, then there is a Type 1 
edge from $Q_{a,b,c}$ to an element below it.  For example, if 
we choose some $d$ such that $a < d \lessdot cb$ and factor
$d$ uniquely as $c'b'$ where $c'\in W^J$ and $u'\in W_J$, then
$Q_{a, b', c'} \lessdot Q_{a,b,c}$.  Therefore by induction,
the lexicographically minimal chain from $Q_{x,u,w}$ to a 
$0$-cell consists entirely of Type 1 edges.  Moreover, it is increasing,
by the argument that we used in the third paragraph of this proof.
Therefore by adding to $m$ the final edge to $\hat{0}$,
we have found an increasing chain from $Q_{x,u,w}$ to $\hat{0}$
which is lexicographically minimal in this interval.

It remains to show that there is a {\it unique} increasing chain
from $Q_{x,u,w}$ to $\hat{0}$.  Suppose that there are two.  
Both of these chains end with the label $\emptyset$, so
by our choice of ordering, these chains must consist of 
Type 1 edges followed by the $\emptyset$ edge.  Therefore
the two increasing chains must both pass through $0$-cells
of the form $Q_{x,b,c}$ and $Q_{x,b',c'}$, where 
$c, c' \in W^J$ and $b, b'\in W_J$.  Since these 
are $0$-cells, we must have  
$x = cb$ and $x = c'b'$.  It follows that $c=c'$ and $b=b'$,
so the increasing chain from $Q_{x,u,w}$ to $\hat{0}$ is indeed 
unique.

This completes the proof that
$\Q^J$ is EL-shellable.
\end{proof}

Recall that $Q_{u_0, u_0, w_0^J}$ is the unique maximal element of 
$\Q^J$.
We now apply Theorem \ref{CW} and 
Proposition \ref{BjornerCW} to the poset 
$P:= \Q^J \setminus \{ Q_{u_0, u_0, w_0^J} \}$.

\begin{corollary}\label{happy}
$\Vert \Q^J \setminus \{ Q_{u_0, u_0, w_0^J} \}\Vert$ is homeomorphic
to a sphere.  Moreover, the 
poset 
$\Q^J \setminus \{ Q_{u_0, u_0, w_0^J} \}$
is the face poset of a regular CW complex which
is homeomorphic to a sphere.
\end{corollary}

Because $\Q^J$ has a unique greatest element, 
Corollary \ref{happy} implies that $\Q^J$ is the face poset of a 
regular CW complex which is homeomorphic to a ball.

Now we apply Lemma \ref{Eulerian} to $\Q^J$.

\begin{corollary} 
$\Q^J$ is Eulerian.
\end{corollary}

Because Rietsch has proved that the closure of every cell in 
$\PP^J_{\geq 0}$ is a union of cells 
(see Theorem \ref{RTheorem}), we can deduce the following.

\begin{corollary}\label{Euler}
The Euler characteristic of the closure of every cell in the 
cell decomposition of $\PP^J_{\geq 0}$ is $1$.
\end{corollary}

These results lead us to make the following conjecture.

\begin{conjecture}
The totally nonnegative part of an arbitrary flag variety 
together with its cell decomposition is a regular
CW complex homeomorphic to a ball.
\end{conjecture}

\appendix
\section{$\Q^J$ for the type A Grassmannian}\label{A}

In independent work, Postnikov \cite{Postnikov} has studied the 
poset of cells of a natural 
cell decomposition of 
the totally nonnegative part of the (type A) Grassmannian $Gr_{k,n}^+$, 
and showed that
this poset can be described in terms of certain tableaux (the 
so-called $\Le$-diagrams) and also 
in terms of certain ``decorated" permutations.  
He defined $Gr_{k,n}^+$ to be the subset of the real Grassmannian where
all Plucker coordinates are non-negative, and defined cells to be
subsets of $Gr_{k,n}^+$ with a given vanishing pattern of Plucker
coordinates.  It is not too hard to see that in the case of the 
Grassmannian, Postnikov's cell decomposition is a special case of Rietsch's 
cell decomposition
\cite{Rietsch1}; although we will not prove that here, we will 
give bijections between Rietsch's cells $Q_{x,u,w}$, $\Le$-diagrams, 
and decorated permutations.
Additionally, we will describe in detail the case of the 
Grassmannian $Gr_{2,4}(\R)$.

\begin{figure}[h]
\centerline{\epsfig{figure=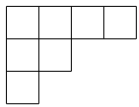}}
\caption{A Young diagram of shape $(4,2,1)$}
\label{YoungDiagram}
\end{figure}

Recall that a {\it partition} $\lambda = (\lambda_1, \dots , \lambda_k)$
is a weakly decreasing sequence of nonnegative numbers.  
For a partition $\lambda$, where $\sum \lambda_i = n$, 
the {\it Young diagram} $Y_{\lambda}$ of 
shape $\lambda$ is 
a left-justified diagram of $n$ boxes, with $\lambda_i$ boxes
in the $i$th row.
Figure \ref{YoungDiagram} shows 
a Young diagram of shape $(4,2,1)$.

\begin{definition}
Fix $k$ and $n$.  A {\it $\Le$-diagram} 
\footnote{The symbol $\Le$ is meant to remind the reader of the 
shape of the forbidden pattern, and should be pronounced as
[{le}], because of its relationship to the letter $L$.}  
$(\lambda, D)_{k,n}$
is a partition $\lambda$ contained in a $k \times (n-k)$ rectangle,
together with a filling $D: Y_{\lambda} \to \{0,+\}$ which has the 
{\it $\Le$-property}:
there is no $0$ which has a $+$ above it and a $+$ to its
left.
\end{definition}

(Here, ``above" means above and in the same column, and 
``to its left" means to the left and in the same row.)
In Figure \ref{LeDiagram} we give 
an example of a $\Le$-diagram.

\begin{figure}[h]
\pspicture(100,0)(230,82)

\rput(190,36)
{$\begin{array}{l}
k=6,\ n=17\\
\lambda=(10,9,9,8,5,2)
\end{array}$}

\rput(-10,36){$k$}
\rput(60,82){$n-k$}
\psline[linecolor=black,linewidth=0.5pt]{-}(0,0)(132,0)(132,72)(0,72)(0,0)
\tbox(0,0){+}
\tbox(1,0){+}

\tbox(0,1){0}
\tbox(1,1){0}
\tbox(2,1){0}
\tbox(3,1){+}
\tbox(4,1){+}

\tbox(0,2){0}
\tbox(1,2){0}
\tbox(2,2){0}
\tbox(3,2){0}
\tbox(4,2){0}
\tbox(5,2){0}
\tbox(6,2){+}
\tbox(7,2){+}

\tbox(0,3){0}
\tbox(1,3){0}
\tbox(2,3){0}
\tbox(3,3){0}
\tbox(4,3){0}
\tbox(5,3){0}
\tbox(6,3){0}
\tbox(7,3){0}
\tbox(8,3){0}

\tbox(0,4){+}
\tbox(1,4){+}
\tbox(2,4){+}
\tbox(3,4){+}
\tbox(4,4){0}
\tbox(5,4){+}
\tbox(6,4){+}
\tbox(7,4){+}
\tbox(8,4){+}

\tbox(0,5){0}
\tbox(1,5){+}
\tbox(2,5){+}
\tbox(3,5){0}
\tbox(4,5){0}
\tbox(5,5){+}
\tbox(6,5){0}
\tbox(7,5){+}
\tbox(8,5){0}
\tbox(9,5){+}

\endpspicture

\caption{A $\Le$-diagram $(\lambda, D)_{k,n}$}
\label{LeDiagram}
\end{figure}

The {\it rank} of 
$(\lambda, D)_{k,n}$ is the number of 
$+$'s in the filling $D$.

\begin{definition}
A {\it decorated permutation} $\tilde{\pi} = (\pi, d)$ is 
a permutation $\pi$ in the symmetric group
$S_n$ together with a coloring (decoration)
$d$ of its fixed points $\pi (i)=i$ by two colors, 
``clockwise" and 
``counterclockwise."
\end{definition}

We represent a decorated permutation $\tilde{\pi}=(\pi,d)$,
where $\pi \in S_n$, by its 
{\it chord diagram}, constructed as follows.  Put $n$ equally
spaced points around a circle, and label these points from $1$ to
$n$ in clockwise order.
If $\pi (i) = j$ then this is represented as a directed arrow, or 
chord,  
from $i$ to $j$.  If $\pi (i) = i$ then
we draw a chord from $i$ to $i$ (i.e. a loop), and orient it either 
clockwise or counterclockwise, according to $d$.

\begin{figure}[h]
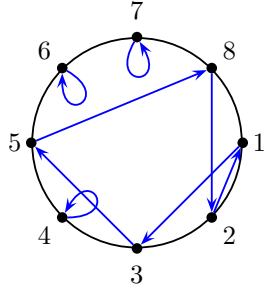

\begin{center}
\pspicture(-60, -60)(60,60)
\pscircle[linecolor=black](0,0){40}
\cnode*[linewidth=0, linecolor=black](40,0){2}{1}
\cnode*[linewidth=0, linecolor=black](28.28,-28.28){2}{2}
\cnode*[linewidth=0, linecolor=black](0, -40){2}{3}
\cnode*[linewidth=0, linecolor=black](-28.28,-28.28){2}{4}
\cnode*[linewidth=0, linecolor=black](-40,0){2}{5}
\cnode*[linewidth=0, linecolor=black](-28.28,28.28){2}{6}
\cnode*[linewidth=0, linecolor=black](0,40){2}{7}
\cnode*[linewidth=0, linecolor=black](28.28,28.28){2}{8}
\rput(46.30,0){$1$}
\rput(35,-35){$2$}
\rput(0,-50){$3$}
\rput(-35,-35 ){$4$}
\rput(-46.30,0){$5$}
\rput(-35,35){$6$}
\rput(0, 50 ){$7$}
\rput(35,35){$8$}
\ncline{->}{1}{3}
\ncline{->}{2}{1}
\ncline{->}{3}{5}
\ncline{->}{5}{8}
\ncline{->}{8}{2}
\nccurve[angleA=-30,angleB=-90,ncurv=20]{->}{6}{6}
\nccurve[angleA=0,angleB=60,ncurv=20]{->}{4}{4}
\nccurve[angleA=-120,angleB=-60,ncurv=20]{->}{7}{7}
\endpspicture
\end{center}
\caption{A chord diagram for a decorated permutation}
\label{chorddiagram}
\end{figure}

For example, the decorated permutation 
$(3,1,5,4,8,6, 7, 2)$ (written in list notation)
with the fixed points $4$, $6$, and $7$ colored in
counterclockwise, clockwise, and counterclockwise, 
respectively, is represented by 
the chord diagram in Figure \ref{chorddiagram}.

Recall that $i$ is a {\it weak excedence} of a permutation $\pi$ 
if $\pi(i) \geq i$.  This definition can be extended to decorated 
permutations as follows: 
$i$ is a {\it weak excedence} of a decorated permutation $(\pi, d)$
if either $\pi(i) > i$ or if $\pi(i)$ is a counterclockwise loop.

We will not review here the rank function on decorated permutations 
nor the order relations for these objects; for details, see 
\cite{Williams}.  However, we should
recall Postnikov's result \cite{Postnikov}
relating the $\Le$-diagrams and decorated permutations to $Gr_{k,n}^+$.

\begin{theorem} \cite{Postnikov}
There is an order-preserving bijection between the poset of cells of 
$Gr_{k,n}^+$ and the poset of $\Le$-diagrams $(\lambda, D)_{k,n}$.
Additionally, there is an order-preserving bijection between the 
poset of cells of $Gr_{k,n}^+$ and the poset of 
decorated permutations on $n$ letters with $k$ weak excedences.
\end{theorem}

We now let 
$W$ be the symmetric group on $n$ letters, 
$S=\{s_1, s_2, \dots , s_{n-1}\}$ be the set of adjacent transpositions, and
$J = \{s_1, s_2, \dots ,\hat{s}_{n-k}, \dots , s_{n-1}\}$.
We will now use $\Q^J$ to denote 
the poset of cells defined in terms of this data.

\begin{lemma}
There is an order-preserving bijection $\Phi_1$ from 
$\Q^J$ to the poset of decorated permutations in $S_n$ with $k$ weak 
excedences, which is defined as follows.
Let $Q_{x,u,w} \in \Q^J$.  Then $\Phi_1(Q_{x,u,w}) = (\pi, d)$ 
where $\pi = xu^{-1}w^{-1}$.  
To define $d$, we make any fixed point that occurs
in one of the positions $w(1), w(2), \dots, w(n-k)$
a {\it clockwise loop}, and we make any fixed point that occurs in one of the
positions 
$w(n-k+1), \dots , w(n)$ a {\it counterclockwise loop}.
\end{lemma}

Additionally, there is a natural bijection between  
$\Le$-diagrams $(\lambda, D)_{k,n}$ and $\Q^J$: 
we thank Postnikov for explaining this to us.

\begin{lemma} \cite{Postnikov2}
There is an order-preserving 
bijection $\Phi_2$  from the
set of $\Le$-diagrams $(\lambda, D)_{k,n}$ to $\Q^J$, defined as follows.
\begin{enumerate}
\item Take $(\lambda,D)_{k,n}$ and replace each $+$
 with an elbow joint $\textelbow$, and each $0$ with a cross
$\textcross$.
\item Note that the west and north borders, and the south and east 
borders, respectively, of $\lambda$, give rise to two length-$n$
paths from the north-east corner to the south-east corner of the 
$k \times (n-k)$ rectangle.  Label each of these paths with the 
numbers $1$ through $n$.
\item View the resulting ``pipe dream" 
as a permutation $w_{\lambda, D} \in S_n$, as in \ref{toy}.
\item Repeat this procedure for the $\Le$-diagram $(\lambda, D_0)_{k,n}$,
where $D_0$ denotes the filling of $\lambda$ which consists 
entirely of $0$'s.  Denote the resulting permutation by 
$w_{\lambda, D_0}$; this permutation is in $W^J$, i.e. it is a Grassmannian 
permutation.
\item Let w:= $w_{\lambda,D_0}$.  Since every element of $W$ has 
  a unique factorization as the product of an element in 
  $W^J_{max}$ and an element in $W_J$, define $x$ and $u$ by 
  the equation $xu^{-1} = w_{\lambda, D}$.  
  We now set $\Phi_2 ((\lambda,D)_{k,n}) := Q_{x,u,w}$.
\end{enumerate}
\end{lemma}

\begin{figure}[h]
\includegraphics[height=.85in]{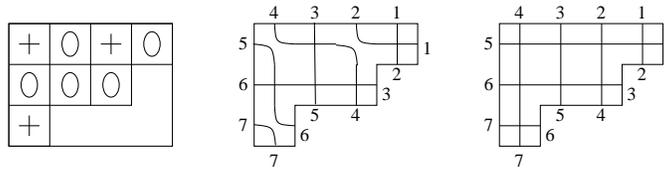}
\caption{The bijection $\Phi_2$}
\label{Bij2}
\end{figure}

\begin{example}\label{toy}
Figure \ref{Bij2}  shows a $\Le$-diagram $(\lambda,D)_{3,7}$
together with the related 
pipe dreams.  This gives rise to the permutation 
$w_{\lambda,D}:= (2,1,5,4,6,3,7)$ (written in list notation) and 
the permutation $w_{\lambda,D_0}:= (2,4,5,7,1,3,6)$.  
\end{example}

For a simple bijection between $\Le$-diagrams and decorated permutations
(which is equal to $\Phi_2 \circ \Phi_1$), 
see \cite{WilliamsStein}.

\begin{example} \label{G24}
We now explain the case of the Grassmannian $Gr_{2,4}(\R)$ in detail.
In that case, the Weyl group $W$ is $S_4$, the symmetric group on 
$4$ letters, and the set $S$ of simple reflections is 
$\{s_1, s_2, s_3\}$ where $s_i$ is the transposition $(i,i+1)$ which
exchanges $i$ and $i+1$.
The subset $J$ is $\{s_1, s_3\}$ and the 
parabolic subgroup $W_J$ is $\langle s_1, s_3 \rangle$.  
A reflection order which puts elements of $W_J$ at the end is 
the following: 
\begin{equation*}
(23) \prec (24) \prec (13) \prec (14) \prec (34) \prec (12)
\end{equation*}

In Figure \ref{Picture}, we have drawn the Hasse diagram of the 
poset $\Q^J$ for the Grassmannian $Gr_{2,4}(\R)$.  Elements $Q_{x,u,w}$ 
(where $x \in W^J_{max}, u\in W_J, w\in W^J$) are
represented by $\Le$-diagrams, and below each $\Le$-diagram, we have 
listed the
triple $(x, u, w)$ corresponding to $Q_{x,u,w}$.  Note that 
in each of these triples we have abbreviated $s_i$ by $i$.
Also note that we have labelled the unique increasing chain from 
the greatest element to the least element of $\Q^J$; every element
in this chain is the totally positive part of a Schubert variety.

\begin{figure}[h]
\includegraphics[height=5.5in, angle=90]{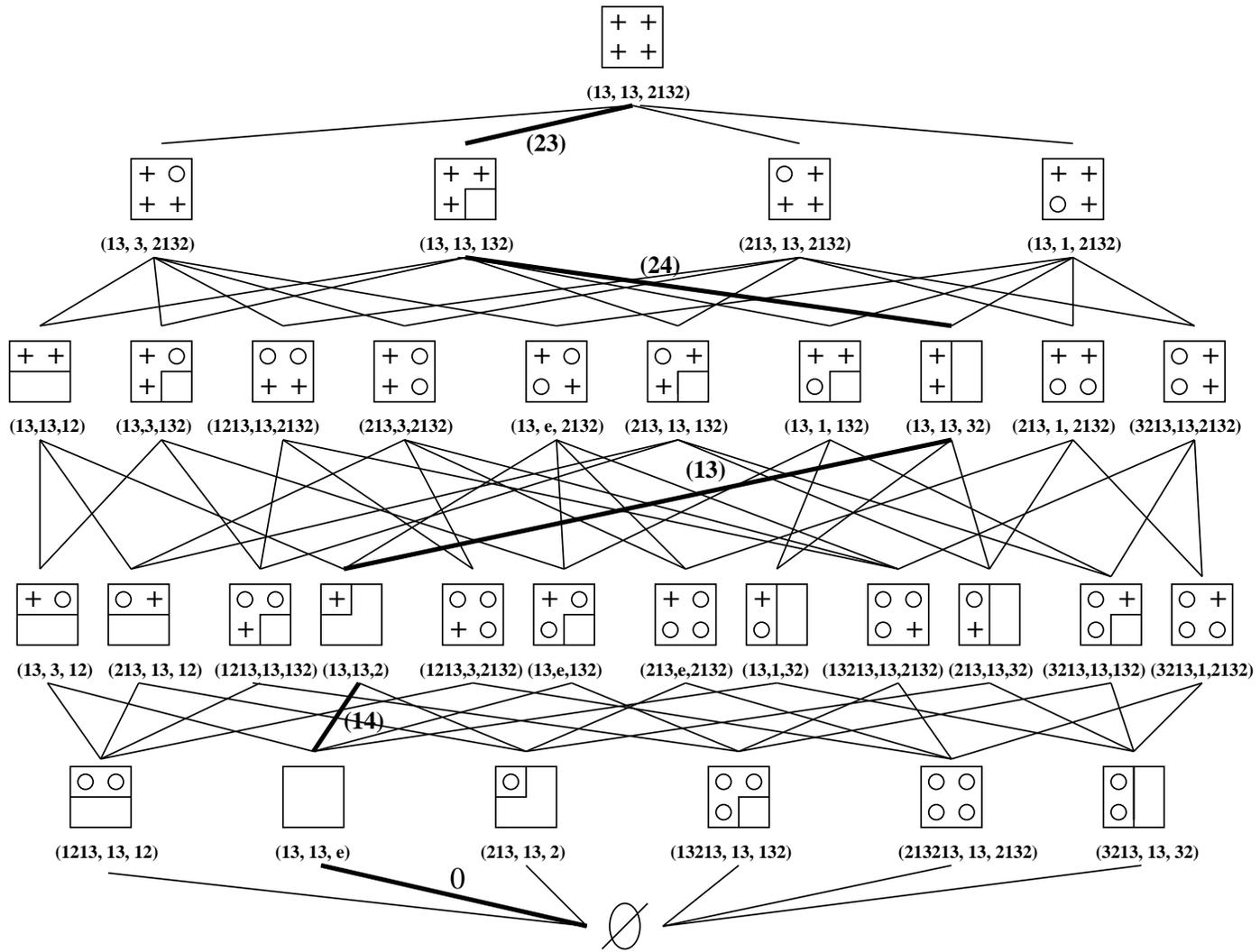}
\caption{$\Q^J$ for the Grassmannian $Gr_{2,4}(\R)$}
\label{Picture}
\end{figure}
\end{example}


\end{document}